\documentclass[12pt,a4paper]{article}
\usepackage{amsmath}
\usepackage{amssymb}
\usepackage{t1enc}
\usepackage[latin1]{inputenc}
\usepackage[english]{babel}
\usepackage{amsfonts}
\usepackage{latexsym}
\usepackage{bm}
\newtheorem{theorem}{Theorem}[section]

\newtheorem{lemma}[theorem]{Lemma}
\newtheorem{cor}[theorem]{Corollary}

\newtheorem{rem}[theorem]{Remark}

\begin{document}
\title{Cross-intersecting integer sequences}

\author{Peter Borg\\[5mm]
Department of Mathematics, University of Malta\\
\texttt{peter.borg@um.edu.mt}}
\date{} 
\maketitle

\begin{abstract}
We call $(a_1, \dots, a_n)$ an \emph{$r$-partial sequence} if exactly $r$ of its entries are positive integers and the rest are all zero. For ${\bf c} = (c_1, \dots, c_n)$ with $1 \leq c_1 \leq \dots \leq c_n$, let $S_{\bf c}^{(r)}$ be the set of $r$-partial sequences $(a_1, \dots, a_n)$ with $0 \leq a_i \leq c_i$ for each $i$ in $\{1, \dots, n\}$, and let $S_{\bf c}^{(r)}(1)$ be the set of members of $S_{\bf c}^{(r)}$ which have $a_1 = 1$. We say that $(a_1, \dots, a_n)$ \emph{meets} $(b_1, \dots, b_m)$ if $a_i = b_i \neq 0$ for some $i$. Two sets $A$ and $B$ of sequences are said to be \emph{cross-intersecting} if each sequence in $A$ meets each sequence in $B$. Let ${\bf d} = (d_1, \dots, d_m)$ with $1 \leq d_1 \leq \dots \leq d_m$. Let $A \subseteq S_{\bf c}^{(r)}$ and $B \subseteq S_{\bf d}^{(s)}$ such that $A$ and $B$ are cross-intersecting. We show that $|A||B| \leq |S_{\bf c}^{(r)}(1)||S_{\bf d}^{(s)}(1)|$ if either $c_1 \geq 3$ and $d_1 \geq 3$ or ${\bf c} = {\bf d}$ and $r = s = n$. 
We also determine the cases of equality.  
We obtain this by proving a general cross-intersection theorem for \emph{weighted} sets. The bound generalises to one for $k \geq 2$ cross-intersecting sets. 
\end{abstract}

\section{Introduction} \label{Def}

Unless otherwise stated, we shall use small letters such as $x$ to
denote elements of a set or non-negative integers or functions,
capital letters such as $X$ to denote sets, and calligraphic
letters such as $\mathcal{F}$ to denote \emph{families}
(i.e.~sets whose elements are sets themselves). It is to be
assumed that arbitrary sets and families are \emph{finite}. We call a set $A$ an \emph{$r$-element set}, or simply an \emph{$r$-set}, if
its size $|A|$ is $r$. For a set $X$, the \emph{power set of $X$} (i.e.~the family of all subsets of $X$) is denoted by $2^X$, and the family of all $r$-element subsets of $X$ is denoted by ${X \choose r}$. 
The set $\{1, 2, \dots\}$ of all positive integers is denoted by $\mathbb{N}$. 
For any $m,n \in \mathbb{N}$ with $m < n$, the set $\{i \in \mathbb{N} \colon m \leq i \leq n\}$ is denoted by $[m,n]$. We abbreviate $[1,n]$ to $[n]$.

We say that a set $A$ \emph{intersects} a set $B$ if $A$ and $B$ contain at least one common element. A family $\mathcal{A}$ of sets is said to be \emph{intersecting} if every two sets in $\mathcal{A}$ intersect. If $x$ is an element of at least one set in a family $\mathcal{F}$, then we call the family of all the sets in $\mathcal{F}$ that contain $x$ the \emph{star of $\mathcal{F}$} with \emph{centre} $x$. A star of a family is the simplest example of an intersecting subfamily.

One of the most popular endeavours in extremal set theory is that of determining the size of a largest intersecting subfamily of a given family $\mathcal{F}$. This took off with \cite{EKR}, which features the classical result, known as the Erd\H os-Ko-Rado (EKR) Theorem, that says that if $r \leq n/2$, then the size of a largest intersecting subfamily of ${[n] \choose r}$ is the size ${n-1 \choose r-1}$ of every star of ${[n] \choose r}$. If $r < n/2$, then, by the Hilton-Milner Theorem \cite{HM}, $\mathcal{A}$ attains the bound if and only if $\mathcal{A}$ is a star of ${[n] \choose r}$. There are various proofs of the EKR Theorem, two of which are particularly short and beautiful: Katona's \cite{K}, introducing the elegant cycle method, and Daykin's \cite{D}, using the fundamental Kruskal-Katona Theorem \cite{Ka,Kr}. Various generalisations and analogues have been obtained. In particular, a sequence of results \cite{EKR,F_t1,W,AK1} culminated in the solution of the more general problem of determining each subfamily $\mathcal{A}$ of ${[n] \choose r}$ whose size is maximum under the condition that every two sets in $\mathcal{A}$ have at least $t$ common elements. The same $t$-intersection problem for $2^{[n]}$ was solved in \cite{Kat}. These are among the most prominent results in extremal set theory. Among the many analogues of the EKR Theorem we find a significant number for families of integer sequences \cite{Kl,B,Meyer,L,FF2,DF,G,E,BL,AK2,FT2,Bey1,Bey2,Borg1,Borg,Borg6}. 
The EKR Theorem inspired a wealth of results that establish how large a system of sets can be under certain intersection conditions; see \cite{DF,F,F2,Borg7}.

Families $\mathcal{A}_1, \dots, \mathcal{A}_k$ are said to be \emph{cross-intersecting} if for every $i$ and $j$ in $[k]$ with $i \neq j$, each set in $\mathcal{A}_i$ intersects each set in $\mathcal{A}_j$.

For intersecting subfamilies of a given family $\mathcal{F}$,
the natural question to ask is how large they can be. For
cross-intersecting families, two natural parameters arise: the
sum and the product of sizes of the cross-intersecting
families (note that the product of sizes of $k$ families
$\mathcal{A}_1, \dots, \mathcal{A}_k$ is the number of $k$-tuples
$(A_1, \dots, A_k)$ such that $A_i \in \mathcal{A}_i$ for each $i
\in [k]$). It is therefore natural to consider the problem of
maximising the sum or the product of sizes of $k$ cross-intersecting subfamilies (not necessarily distinct or non-empty) of a given family $\mathcal{F}$. In \cite{Borg8} this problem is analysed in general, and it is shown that for $k$ sufficiently large it reduces to the intersection problem (i.e.~the problem of maximising the size of an intersecting subfamily of $\mathcal{F}$). 
Solutions have been obtained for various families (see \cite{Borg8}), including ${[n] \choose r}$ \cite{H,Pyber, MT,Bey3,Borg4,WZ}, $2^{[n]}$ \cite{MT2,Borg8} and families of integer sequences \cite{Moon,Borg,WZ,Zhang,Tok3}. It is worth pointing out that Wang and Zhang \cite{WZ} solved the maximum sum problem for an important class of families that includes ${[n] \choose r}$ and families of integer sequences, using a striking combination of the method in \cite{Borg4,Borg3,Borg2,BL2,Borg5} and an important lemma that is found in \cite{AC,CK} and referred to as the `no-homomorphism lemma'. The solution for ${[n] \choose r}$ had been obtained by Hilton \cite{H} and is the first result of this kind. In this paper we are concerned with the maximum product problem for integer sequences. We will now introduce a more general problem. 

We will represent a sequence $a_1, \dots, a_n$ by an $n$-tuple $(a_1, \dots, a_n)$, and we say that it is of \emph{length $n$}. In this paper all sequences are taken to be sequences of non-negative integers. We call a sequence of positive integers a \emph{positive sequence}. We call $(a_1, \dots, a_n)$ an \emph{$r$-partial sequence} if exactly $r$ of its entries are positive integers and the rest are all zero, i.e.~if $|\{i \in [n] \colon a_i \neq 0\}| = r$. Thus, an $n$-partial sequence of length $n$ is positive. A sequence $(c_1, \dots, c_n)$ is said to be \emph{increasing} if $c_1 \leq \dots \leq c_n$. We call an increasing positive sequence an \emph{IP sequence}. Note that $(c_1, \dots, c_n)$ is an IP sequence if and only if $1 \leq c_1 \leq \dots \leq c_n$. 

For an IP sequence ${\bf c} = (c_1, \dots, c_n)$, let 
\[S_{\bf c}^{(r)} = \{(a_1, \dots, a_n) \colon (a_1, \dots, a_n) \mbox{ is an $r$-partial sequence}, \, 0 \leq a_i \leq c_i \mbox{ for each } i \in [n]\}\] 
and let $S_{\bf c}^{(r)}(1)$ denote the set $\{(a_1, \dots, a_n) \in S_{\bf c}^{(r)} \colon a_1 = 1\}$. Note that $S_{\bf c}^{(n)}$ is the Cartesian product $[c_1] \times \dots \times [c_n]$, and we may abbreviate this to $S_{\bf c}$. 

We say that ${\bf a} = (a_1, \dots, a_n)$ \emph{meets} ${\bf b} = (b_1, \dots, b_m)$ if the two sequences agree on some positive entry, i.e.~if $a_i = b_i \neq 0$ for some $i$. If ${\bf a} \in S_{\bf c}$, then $\bf a$ meets $\bf b$ if and only if $a_i = b_i$ for some $i$.

The following is our first result, which we will prove in Sections~\ref{Weightedsection} and \ref{Proofmain}.

\begin{theorem} \label{main} If $\bf{c}$ is an IP sequence $(c_1, \dots, c_n)$ and $A, B \subseteq S_{\bf c}$ such that each sequence in $A$ meets each sequence in $B$, then
\[|A||B| \leq \left(\frac{|S_{\bf c}|}{c_1} \right)^2 = \left(\frac{1}{c_1}\prod_{i=1}^n c_i\right)^2,\]
and equality holds if $A = B = S_{\bf c}(1)$. Moreover, unless $c_1 = 2$, the bound is attained if and only if $A = B = \{(a_1, \dots, a_n) \in S_{\bf c} \colon a_p = q\}$ for some $p \in \{i \in [n] \colon c_i = c_1\}$ and some $q \in [c_p]$.
\end{theorem}
For $c_1 = 2$ there may be other optimal configurations; for example, if $c_1 = c_2 = c_3 = 2$ and $A = B = \{(a_1, \dots, a_n) \in S_{\bf c} \colon |\{i \in [3] \colon a_i = 1\}| \geq 2\}$, then $|A||B|$ is maximum too. 

The EKR-type version of Theorem~\ref{main} is the solution to the problem of maximising the size of a subset $A$ of $S_{\bf c}$ under the condition that every two sequences in $A$ meet, and this is given in \cite{B,L,Borg} (for $c_1 = c_n$ this is given in a stronger form in \cite{FF2,AK2,FT2}); 
this follows from Theorem~\ref{main} by taking $A = B$. The special case $c_1 = c_n$ of Theorem~\ref{main} has already been treated by Moon \cite{Moon}, Tokushige \cite{Tok3} and Zhang \cite{Zhang}. Moon proved it in a stronger form for $c_1 \geq 3$, using a non-trivial induction argument. Tokushige proved it for $c_1 \geq 4$, also in a stronger form, using an eigenvalue method. Zhang proved it for $c_1 \geq 4$, using Katona's cycle method. Allowing $\bf c$ to be increasing appears to be a significant relaxation for the product problem, and in fact one can see from Zhang's proof why the argument there does not carry forward to this more general setting. Our approach will be based on the idea of generalising the setting enough for induction to work, and we will use the \emph{compression} technique in two different ways (see Sections~\ref{Compsection} and \ref{Proofmain}), together with a new alteration method. For this to work for the case $c_1 \geq 3$, it suffices to consider the generalisation that $A \subseteq S_{\bf c}$ and $B \subseteq S_{\bf d}$, where $\bf c$ and $\bf d$ are arbitrary IP sequences whose first entry is at least 3. However, our method allows us to prove the following generalisation for $r$-partial sequences (the proof of which is also given in Sections~\ref{Weightedsection} and \ref{Proofmain}). 

\begin{theorem} \label{mainpartial} Let ${\bf c} = (c_1, \dots, c_m)$ and ${\bf d} = (d_1, \dots, d_n)$ be IP sequences such that $c_1 \geq 3$ and $d_1 \geq 3$. Let $r \in [m]$ and $s \in [n]$. If $A \subseteq S_{\bf c}^{(r)}$, $B \subseteq S_{\bf d}^{(s)}$, and each sequence in $A$ meets each sequence in $B$, then
\[|A||B| \leq \Bigg{(} \sum_{I \in {[2,m] \choose r-1}} \prod_{i \in I} c_i \Bigg{)} \Bigg{(} \sum_{J \in {[2,n] \choose s-1}}\prod_{j \in J} c_j \Bigg{)},\]
and equality holds if and only if, for some $p \in \{h \in [\min\{m,n\}] \colon c_h = c_1, d_h = d_1\}$ and some $q \in [c_p]$, $A = \{(a_1, \dots, a_m) \in S_{\bf c}^{(r)} \colon a_p = q\}$ and $B = \{(b_1, \dots, b_n) \in S_{\bf d}^{(s)} \colon b_p = q\}$.
\end{theorem}
Thus, for $c_1 \geq 3$, Theorem~\ref{main} is the special case when $r = s = m = n$ and ${\bf c} = {\bf d}$. In general, Theorem~\ref{mainpartial} does not hold for $c_1 = 1$; indeed, if $c_1 = c_m = d_1 = d_n = 1$, $m = n$ and $m/2 < r = s < m$, then any two sequences in $S_{\bf c}^{(r)}$ meet and hence we can take $A = B = S_{\bf c}^{(r)}$. The case when $c_1 = 2$ or $d_1 = 2$ seems to require special treatment and remains a problem to be investigated.

The EKR-type version of Theorem~\ref{mainpartial} is the solution to the problem of maximising the size of a subset $A$ of $S_{\bf c}^{(r)}$ under the condition that every two sequences in $A$ meet, and this is given in \cite{DF,E,BL,Bey2} 
for $c_1 = c_n \geq 2$, in \cite{BE,HST} for $c_1 \geq 2$, and in \cite{Bey2} for $c_1 = 1$; for $c_1 \geq 3$ this follows from Theorem~\ref{mainpartial} by taking ${\bf c} = {\bf d}$, $r = s$ and $A = B$.

Theorems~\ref{main} and \ref{mainpartial} are consequences of a result in Section~\ref{Weightedsection} for cross-intersecting families of \emph{weighted} sets, the proof of which contains the main ideas and new observations in this paper. Basically, the method is as follows. We consider two cross-intersecting families $\mathcal{A}$ and $\mathcal{B}$, where $\mathcal{A}$ consists of weighted subsets of $[m]$, $\mathcal{B}$ consists of weighted subsets of $[n]$, and some conditions hold. We use induction on $m$ and $n$. The challenging part is the case $m = n$. The first problem that arises is that we can have a set $A \in \mathcal{A}$ and a set $B \in \mathcal{B}$ that intersect only in $n$; in this case, we cannot simply remove $n$ and apply the induction hypothesis. Thus, we consider two alterations: removing $A$ from $\mathcal{A}$ and adding $B \backslash \{n\}$ to $\mathcal{B}$, and removing $B$ from $\mathcal{B}$ and adding $A \backslash \{n\}$ to $\mathcal{A}$. This yields two new pairs of cross-intersecting families. The second problem is that the product of the weights of a new pair obtained in this way may become smaller. The critical part of the proof is the observation that if we assume that this happens for both pairs, then remarkably the resulting inequalities lead to a contradiction. Thus, we can replace the original pair by the new pair that gives a product that is not smaller than the original one. By repeated application of this alteration, we eventually obtain a pair of cross-intersecting families for which the first problem does not arise. For convenience, we will actually tackle the two problems in a different but equivalent way. We will take $\mathcal{A}$ and $\mathcal{B}$ to be such that the product of their weights is maximum. Then we show that if we assume that the first problem arises, then, by considering the alterations above, we can construct two new cross-intersecting families for which the product of weights is larger than that for $\mathcal{A}$ and $\mathcal{B}$, a contradiction.

Each of Theorems~\ref{main} and \ref{mainpartial} generalises to one for any $k \geq 2$ subsets. Together they generalise as follows.

\begin{theorem} \label{mainpartialgen} Let ${\bf c}_1 = (c_{1,1}, \dots, c_{1,n_1}), \dots, {\bf c}_k = (c_{k,1}, \dots, c_{k,n_k})$ be IP sequences. For each $i \in [k]$, let $r_i \in [n_i]$. Suppose that either $c_{i,1} \geq 3$ for each $i \in [k]$ or ${\bf c}_1 = \dots = {\bf c}_k$ and $r_1 = \dots = r_k = n_1$. Let $A_1 \subseteq S_{{\bf c}_1}^{(r_1)}, \dots, A_k \subseteq S_{{\bf c}_k}^{(r_k)}$ such that for every $p, q \in [k]$ with $p \neq q$, each sequence in $A_p$ meets each sequence in $A_q$. Then
\[\prod_{i = 1}^k |A_i| \leq \prod_{i=1}^k \Bigg{(} \sum_{I \in [2,n_i] \choose r_i-1}\prod_{j \in I} c_{i,j} \Bigg{)},\]
and equality holds if $A_i = S_{{\bf c}_i}^{(r_i)}(1)$ for each $i \in [k]$. Moreover, unless $c_{1,1} = 2$, the bound is attained if and only if, for some $p \in \{h \in [\min\{n_1, \dots, n_k\}] \colon c_{i,h} = c_{i,1} \mbox{ for each } i \in [k]\}$ and some $q \in [c_p]$, $A_i = \{(a_1, \dots, a_{n_i}) \in S_{{\bf c}_i}^{(r_i)} \colon a_p = q\}$ for each $i \in [k]$.
\end{theorem}
\textbf{Proof.} For each $i \in [k]$, let $x_i = |A_i|$ and $y_i = |S_{{\bf c}_i}^{(r_i)}(1)|$. By Theorems~\ref{main} and \ref{mainpartial}, $x_ix_j \leq y_iy_j$ for any $i, j \in [k]$ with $i \neq j$. Let mod$^*$ be the usual \emph{modulo operation} with the exception that for every two integers $s$ and $t > 0$, $(st) \, {\rm mod}^* \, t$ is $t$ instead of $0$. We have
\begin{align} \left( \prod_{i=1}^k x_i \right)^2 &=
(x_1x_2)(x_{3 \, {\rm mod}^* \, k}x_{4 \, {\rm mod}^* \,
k})\cdots(x_{(2k-1) \, {\rm mod}^* \, k} x_{(2k) \, {\rm mod}^*
\, k}) \nonumber \\
&\leq (y_1y_2)(y_{3 \, {\rm mod}^* \, k}y_{4 \, {\rm mod}^* \,
k})\cdots(y_{(2k-1) \, {\rm mod}^* \, k} y_{(2k) \, {\rm mod}^*
\, k}) = \left( \prod_{i=1}^k y_i \right)^2. \nonumber
\end{align}
So $\prod_{i=1}^k x_i \leq \prod_{i=1}^k y_i$, and this is the bound in the theorem. Suppose equality holds and $c_{1,1} \neq 2$. Then $x_1x_2 = y_1y_2$. By Theorems~\ref{main} and \ref{mainpartial}, $A_1 = \{(a_1, \dots, a_{n_1}) \in S_{{\bf c}_1}^{(r_1)} \colon a_p = q\}$ for some $p \in [n_1]$ and some $q \in [c_p]$. 

Suppose $c_{1,1} = 1$. Then ${\bf c}_1 = \dots = {\bf c}_k$ and $r_1 = \dots = r_k = n_1$. So $S_{{\bf c}_i}^{(r_i)} = S_{{\bf c}_1}^{(n_1)} = S_{{\bf c}_1} = S_{{\bf c}_1}(1)$ for each $i \in [k]$. Thus, the bound in the theorem is attained only if $A_i = S_{{\bf c}_i}^{(r_i)}$ for each $i \in [k]$. 

Now suppose $c_{1,1} \geq 3$. Consider $j \in [2,k]$. Clearly, for each ${\bf a} \in S_{{\bf c}_j}^{(r_j)}$ such that the $p$'th entry of $\bf a$ is not $q$, there are sequences in $A_1$ that do not meet $\bf a$. Taking $B_j = \{(a_1, \dots, a_{n_j}) \in S_{{\bf c}_j}^{(r_j)} \colon a_p = q\}$, we therefore have that $A_j \subseteq B_j$. Clearly, the bound in the theorem is attained only if $A_j = B_j$ for each $j \in [2,k]$ and $c_{i,p} = c_{i,1}$ for each $i \in [k]$.~\hfill{$\Box$} \\

The above results can be phrased in terms of cross-intersecting families of sets as follows.
For any IP sequence ${\bf c} = (c_1, \dots, c_n)$ and any integer $r$, let 
\[\mathcal{L}_{\bf c}^{(r)} = \left\{\{(i_1,a_{i_1}), \dots, (i_r, a_{i_r})\} \colon \{i_1, \dots, i_r\} \in {[n] \choose r}, \, a_{i_j} \in [c_{i_j}] \mbox{ for each } j \in [r] \right \}.\]
We call a set in $\mathcal{L}_{\bf c}^{(r)}$ a \emph{labeled} set (following \cite{Borg}). We may abbreviate $\mathcal{L}_{\bf c}^{(n)}$ to $\mathcal{L}_{\bf c}$; note that $\mathcal{L}_{\bf c} = \{\{(1,a_1), \dots, (n,a_n)\} \colon a_i \in [c_i] \mbox{ for each } i \in [n]\}$.

\begin{rem}\label{seqlab} \emph{There is an obvious one-to-one correspondence between $S_{\bf c}^{(r)}$ and $\mathcal{L}_{\bf c}^{(r)}$. Indeed, let $f \colon S_{\bf c}^{(r)} \rightarrow \mathcal{L}_{\bf c}^{(r)}$ be the function that maps the sequence ${\bf a} = (a_1, \dots, a_n)$ in $S_{\bf c}^{(r)}$ to $L_{\bf a} = \{(i,a_i) \colon i \in [n], \, a_i \neq 0\}$; then $f$ is a bijection. Moreover, a sequence ${\bf a}$ in $S_{\bf c}^{(r)}$ meets a sequence ${\bf b}$ in $S_{\bf d}^{(s)}$ if and only if the corresponding sets $L_{\bf a}$ and $L_{\bf b}$ intersect.}
\end{rem}

Therefore, Theorem~\ref{mainpartialgen} can be re-phrased as follows.

\begin{theorem}[Theorem~\ref{mainpartialgen} re-phrased] \label{setversion} Let ${\bf c}_1, \dots, {\bf c}_k$ and $r_1, \dots, r_k$ be as in Theorem~\ref{mainpartialgen}. Let $\mathcal{A}_1 \subseteq \mathcal{L}_{{\bf c}_1}^{(r_1)}, \dots, \mathcal{A}_k \subseteq \mathcal{L}_{{\bf c}_k}^{(r_k)}$ such that $\mathcal{A}_1, \dots, \mathcal{A}_k$ are cross-intersecting. Then
\[\prod_{i = 1}^k |\mathcal{A}_i| \leq \prod_{i=1}^k \Bigg{(}\sum_{I \in {[2,n_i] \choose r_i-1}}\prod_{j \in I} c_{i,j} \Bigg{)},\]
and equality holds if $\mathcal{A}_i = \{A \in \mathcal{L}_{{\bf c}_i}^{(r_i)} \colon (1,1) \in A\}$ for each $i \in [k]$. Moreover, unless $c_{1,1} = 2$, the bound is attained if and only if, for some $p \in \{h \in [\min\{n_1, \dots, n_k\}] \colon c_{i,h} = c_{i,1} \mbox{ for each } i \in [k]\}$ and some $q \in [c_p]$, $\mathcal{A}_i = \{A \in \mathcal{L}_{{\bf c}_i}^{(r_i)} \colon (p,q) \in A\}$ for each $i \in [k]$.
\end{theorem}
%
For the special case $\mathcal{L}_{{\bf c}_1}^{(r_1)} = \dots = \mathcal{L}_{{\bf c}_k}^{(r_k)} = \mathcal{L}_{\bf c}$, the analogous result for the maximum sum of sizes 
is given in \cite{Borg}. For the case $\mathcal{L}_{{\bf c}_1}^{(r_1)} = \dots = \mathcal{L}_{{\bf c}_k}^{(r_k)} = \mathcal{L}_{\bf c}^{(r)}$ with ${\bf c} = (c, \dots, c)$, the maximum sum of sizes is determined in \cite{BL2}.

We now start working towards the proofs of Theorems~\ref{main} and \ref{mainpartial}. The next section is dedicated to some basic results we need about the compression operation, which is a very useful tool in extremal set theory. Then in Section~\ref{Weightedsection} we prove the cross-intersection result for weighted sets using compressions, and from this we obtain the proofs of Theorems~\ref{main} and \ref{mainpartial}, finalised in Section~\ref{Proofmain}.


\section{The compression operation}
\label{Compsection}
For any $i, j \in [n]$, let $\delta_{i,j}
\colon 2^{[n]} \rightarrow 2^{[n]}$ be defined by
\[ \delta_{i,j}(A) = \left\{ \begin{array}{ll}
(A \backslash \{j\}) \cup \{i\} & \mbox{if $j \in A$ and $i \notin
A$};\\
A & \mbox{otherwise,}
\end{array} \right. \]
and let $\Delta_{i,j} \colon 2^{2^{[n]}} \rightarrow 2^{2^{[n]}}$ be the
\emph{compression operation} (see \cite{EKR}) defined by
\[\Delta_{i,j}(\mathcal{A}) = \{\delta_{i,j}(A) \colon A \in
\mathcal{A}, \delta_{i,j}(A) \notin \mathcal{A}\} \cup \{A \in
\mathcal{A} \colon \delta_{i,j}(A) \in \mathcal{A}\}.\]
Note that $|\Delta_{i,j}(\mathcal{A})| = |\mathcal{A}|$. \cite{F}
provides a survey on the properties and uses of compression (also
called \emph{shifting}) operations in extremal set theory. We will need the following basic result, which we prove for completeness.

\begin{lemma}\label{compcross} Let $\mathcal{A}$ and $\mathcal{B}$ be cross-intersecting subfamilies of $2^{[n]}$, and let $i, j \in [n]$. Then $\Delta_{i,j}(\mathcal{A})$ and $\Delta_{i,j}(\mathcal{B})$ are cross-intersecting subfamilies of $2^{[n]}$.
\end{lemma}
\textbf{Proof.} 
Suppose $A \in \Delta_{i,j}(\mathcal{A})$ and $B \in \Delta_{i,j}(\mathcal{B})$. If $A \in \mathcal{A}$ and $B \in \mathcal{B}$, then $A \cap B \neq \emptyset$ since $\mathcal{A}$ and $\mathcal{B}$ are cross-intersecting. Suppose $A \notin \mathcal{A}$ or $B \notin \mathcal{B}$; we may assume that $A \notin \mathcal{A}$. Then $A = \delta_{i,j}(A') \neq A'$ for some $A' \in \mathcal{A}$.  So $i \notin A'$, $j \in A'$, $i \in A$ and $j \notin A$. Suppose $A \cap B = \emptyset$. So $i \notin B$ and hence $B \in \mathcal{B}$. So $B \in \mathcal{B} \cap \Delta_{i,j}(\mathcal{B})$ and hence $B, \delta_{i,j}(B) \in \mathcal{B}$. So $A' \cap B \neq \emptyset$ and $A' \cap \delta_{i,j}(B) \neq \emptyset$. From $A \cap B = \emptyset$ and $A' \cap B \neq \emptyset$ we get $A' \cap B = \{j\}$, but this yields the contradiction that $A' \cap \delta_{i,j}(B) = \emptyset$.~\hfill{$\Box$}\\


If $i < j$, then we call $\Delta_{i,j}$ a \emph{left-compression}. A family $\mathcal{F} \subseteq 2^{[n]}$ is said to be
\emph{compressed} if $\Delta_{i,j}(\mathcal{F}) = \mathcal{F}$ for any $i,j \in [n]$ with $i < j$. In other words, $\mathcal{F}$ is compressed if it is invariant under left-compressions. Note that $\mathcal{F}$ is compressed if and only if $(F \backslash \{j\}) \cup \{i\} \in \mathcal{F}$ whenever $1 \leq i < j \in F \in \mathcal{F}$ and $i
\in [n] \backslash F$.

Suppose that a subfamily $\mathcal{A}$ of $2^{[n]}$ is not compressed. Then $\mathcal{A}$ can be transformed to a compressed family through left-compressions as follows. Since $\mathcal{A}$ is not compressed, we can find a left-compression that changes $\mathcal{A}$, and we apply it to $\mathcal{A}$ to obtain a new subfamily of $2^{[n]}$. We keep on repeating this (always applying a left-compression to the last family obtained) until we obtain a subfamily of $2^{[n]}$ that is invariant under any left-compression (such a point is indeed reached, because if $\Delta_{i,j}(\mathcal{F}) \neq \mathcal{F} \subseteq 2^{[n]}$ and $i < j$, then $0 < \sum_{G \in \Delta_{i,j}(\mathcal{F})} \sum_{b \in G} b < \sum_{F \in \mathcal{F}} \sum_{a \in F} a$).


Now consider $\mathcal{A}, \mathcal{B} \subseteq 2^{[n]}$ such that $\mathcal{A}$ and $\mathcal{B}$ are cross-intersecting. Then, by Lemma~\ref{compcross}, we can obtain $\mathcal{A}^*, \mathcal{B}^* \subseteq 2^{[n]}$ such that $\mathcal{A}^*$ and $\mathcal{B}^*$ are compressed and cross-intersecting, $|\mathcal{A}^*| = |\mathcal{A}|$ and $|\mathcal{B}^*| = |\mathcal{B}|$. Indeed, similarly to the above procedure, if we can find a left-compression that changes at least one of $\mathcal{A}$ and $\mathcal{B}$, then we apply it to both $\mathcal{A}$ and $\mathcal{B}$, and we keep on repeating this (always performing this on the last two families obtained) until we obtain $\mathcal{A}^*, \mathcal{B}^* \subseteq 2^{[n]}$ such that both $\mathcal{A}^*$ and $\mathcal{B}^*$ are invariant under any left-compression.

\section{A cross-intersection theorem for weighted sets} \label{Weightedsection}

Let $\mathbb{R}^+$ denote the set of positive real numbers.  For any non-empty family $\mathcal{F}$, any function $w \colon \mathcal{F} \rightarrow \mathbb{R}^+$ (which we call a \emph{weight function}), and any $\mathcal{A} \subseteq \mathcal{F}$, we denote
the sum $\sum_{A \in \mathcal{A}} w(A)$ (of \emph{weights} of
sets in $\mathcal{A}$) by $w^{(\mathcal{F})}(\mathcal{A})$. Note that if $\mathcal{A}$ is empty, then $w^{(\mathcal{F})}(\mathcal{A})$ is the \emph{empty sum} and we will adopt the convention of taking this to be $0$.

If $x$ is an element of a set $X$ and $\mathcal{F} \subseteq 2^X$, then we denote the family $\{F \in \mathcal{F} \colon x \in F\}$ by $\mathcal{F}(x)$.

A family $\mathcal{H}$ is said to be \emph{hereditary} if for each $H \in \mathcal{H}$, all the subsets of $H$ are in $\mathcal{H}$. Thus,
a family is hereditary if and only if it is a union of power sets.

\begin{theorem}\label{xintweight} Let $m, n \in \mathbb{N}$. Let $\emptyset \neq \mathcal{G} \subseteq 2^{[m]}$ and $\emptyset \neq \mathcal{H} \subseteq 2^{[n]}$ such that $\mathcal{G}$ and $\mathcal{H}$ are hereditary and compressed. For each $\mathcal{F} \in \{\mathcal{G}, \mathcal{H}\}$, let $w_{\mathcal{F}} \colon \mathcal{F} \rightarrow \mathbb{R}^+$ such that \\
(a) $w_{\mathcal{F}}(A) \geq 2w_{\mathcal{F}}(B)$ for every $A, B \in \mathcal{F}$ with $A \subsetneq B$, and \\
(b) $w_{\mathcal{F}}(\delta_{i,j}(C)) \geq w_{\mathcal{F}}(C)$ for every $C \in \mathcal{F}$ and every $i,j \in \left[ \max\{m,n\} \right]$ with $i < j$. \\
Let $g = w_{\mathcal{G}}$ and $h = w_{\mathcal{H}}$. If $\mathcal{A} \subseteq \mathcal{G}$ and $\mathcal{B} \subseteq  \mathcal{H}$ such that $\mathcal{A}$ and $\mathcal{B}$ are cross-intersecting, then
\[g^{(\mathcal{G})}(\mathcal{A}) h^{(\mathcal{H})}(\mathcal{B}) \leq g^{(\mathcal{G})}(\mathcal{G}(1)) h^{(\mathcal{H})}(\mathcal{H}(1)).\]
Moreover, equality holds if and only if $\mathcal{A} = \mathcal{G}(a)$ and $\mathcal{B} = \mathcal{H}(a)$ for some $a \in [m] \cap [n]$ such that $g^{(\mathcal{G})}(\mathcal{G}(a)) = g^{(\mathcal{G})}(\mathcal{G}(1))$ and $h^{(\mathcal{H})}(\mathcal{H}(a)) = h^{(\mathcal{H})}(\mathcal{H}(1))$.  
\end{theorem}
%
%
%

For the extremal cases, we shall use the following lemma.

\begin{lemma}\label{complemma1} Let $\mathcal{H}$ be a compressed subfamily of $2^{[n]}$. Let $w \colon \mathcal{H} \rightarrow \mathbb{R}^+$ such that  $w(\delta_{i,j}(H)) \geq w(H)$ for every $H \in \mathcal{H}$ and every $i,j \in [n]$ with $i < j$. Then $w^{(\mathcal{H})}(\mathcal{H}(a)) \leq w^{(\mathcal{H})}(\mathcal{H}(1))$ for each $a \in [n]$.
\end{lemma}
\textbf{Proof.} Let $a \in [n]$. Let $\mathcal{D} = \Delta_{1,a}(\mathcal{H}(a))$. Since $\mathcal{H}$ is compressed, $\mathcal{D} \subseteq \mathcal{H}$. Thus it is immediate from the definitions of $\mathcal{D}$ and $w$ that $w^{(\mathcal{H})}(\mathcal{D}) \geq w^{(\mathcal{H})}(\mathcal{H}(a))$. The result follows if we show that $1 \in D$ for each $D \in \mathcal{D}$. Let $D \in \mathcal{D}$. If $D \notin \mathcal{H}(a)$, then $D = \delta_{1,a}(H) \neq H$ for some $H \in \mathcal{H}(a)$, and hence $1 \in D$. Suppose that $D \in \mathcal{H}(a)$. Then, since $D \in \Delta_{1,a}(\mathcal{H}(a))$, $\delta_{1,a}(D) \in \mathcal{H}(a)$. So $a \in \delta_{1,a}(D)$. Since $a \in D$, it follows that $1 \in D$.~\hfill{$\Box$}\\
\\
\textbf{Proof of Theorem~\ref{xintweight}.} 
%
By induction on $m + n$. The basis is $m + n = 2$ with $m = n = 1$, in which case the result is trivial. Now consider $m + n > 2$. We may assume that $m \leq n$. If $m = 1$, then the result is trivial too, so we consider $m \geq 2$. If at least one of $\mathcal{G}$ and $\mathcal{H}$ is $\{\emptyset\}$, then we trivially have $g^{(\mathcal{G})}(\mathcal{A})h^{(\mathcal{H})}(\mathcal{B}) = 0 = g^{(\mathcal{G})}(\mathcal{G}(1))h^{(\mathcal{H})}(\mathcal{H}(1))$. Thus, we will assume that $\mathcal{G} \neq \{\emptyset\}$ and $\mathcal{H} \neq \{\emptyset\}$, meaning that each of $\mathcal{G}$ and $\mathcal{H}$ contain at least one non-empty set. Since $\mathcal{G}$ and $\mathcal{H}$ are hereditary and compressed, we clearly have $\{1\} \in \mathcal{G}$ and $\{1\} \in \mathcal{H}$. So $g^{(\mathcal{G})}(\mathcal{G}(1)) > 0$ and $h^{(\mathcal{H})}(\mathcal{H}(1)) > 0$. Let $\mathcal{A} \subseteq \mathcal{G}$ and $\mathcal{B} \subseteq  \mathcal{H}$ such that $g^{(\mathcal{G})}(\mathcal{A}) h^{(\mathcal{H})}(\mathcal{B})$ is maximum under the condition that $\mathcal{A}$ and $\mathcal{B}$ are cross-intersecting. Since $\mathcal{G}(1)$ and $\mathcal{H}(1)$ are cross-intersecting, it follows that 
\begin{equation} g^{(\mathcal{G})}(\mathcal{A}) h^{(\mathcal{H})}(\mathcal{B}) \geq g^{(\mathcal{G})}(\mathcal{G}(1)) h^{(\mathcal{H})}(\mathcal{H}(1)) > 0. \label{main0.1}
\end{equation} 

We will first show that we may assume that $\mathcal{A}$ and $\mathcal{B}$ are compressed.

As is explained in Section~\ref{Compsection}, we apply left-compressions to $\mathcal{A}$ and $\mathcal{B}$ simultaneously until we obtain two compressed cross-intersecting families $\mathcal{A}^*$ and $\mathcal{B}^*$ such that $|\mathcal{A}^*| = |\mathcal{A}|$ and $|\mathcal{B}^*| = |\mathcal{B}|$. Since $\mathcal{G}$ and $\mathcal{H}$ are compressed, $\mathcal{A}^* \subset \mathcal{G}$ and $\mathcal{B}^* \subset \mathcal{H}$. From (b) we obtain $g^{(\mathcal{G})}(\mathcal{A}) \leq g^{(\mathcal{G})}(\mathcal{A}^*)$ and $h^{(\mathcal{H})}(\mathcal{B}) \leq h^{(\mathcal{H})}(\mathcal{B}^*)$. By the choice of $\mathcal{A}$ and $\mathcal{B}$, we actually have $g^{(\mathcal{G})}(\mathcal{A}) = g^{(\mathcal{G})}(\mathcal{A}^*)$ and $h^{(\mathcal{H})}(\mathcal{B}) = h^{(\mathcal{H})}(\mathcal{B}^*)$. 

Suppose that $\mathcal{A}^* = \mathcal{G}(c)$ and $\mathcal{B}^* = \mathcal{H}(c)$ for some $c \in [m] \cap [n]$ such that $g^{(\mathcal{G})}(\mathcal{G}(c)) = g^{(\mathcal{G})}(\mathcal{G}(1))$ and $h^{(\mathcal{H})}(\mathcal{H}(c)) = h^{(\mathcal{H})}(\mathcal{H}(1))$. Then $g^{(\mathcal{G})}(\mathcal{G}(c)) > 0$ and $h^{(\mathcal{H})}(\mathcal{H}(c)) > 0$. So  
$\mathcal{G}(c) \neq \emptyset$  and $\mathcal{H}(c) \neq \emptyset$. Thus, since $\mathcal{G}$ and $\mathcal{H}$ are hereditary, $\{c\} \in \mathcal{A}^*$ and $\{c\} \in \mathcal{B}^*$. So $\{a\} \in \mathcal{A}$ for some $a \in [m]$, and $\{b\} \in \mathcal{B}$ for some $b \in [n]$. Since $\mathcal{A}$ and $\mathcal{B}$ are cross-intersecting, we have $a = b$, $\mathcal{A} \subseteq \mathcal{G}(a)$ and $\mathcal{B} \subseteq \mathcal{H}(a)$. Since $\mathcal{G}(a)$ and $\mathcal{H}(a)$ are cross-intersecting, it follows by the choice of $\mathcal{A}$ and $\mathcal{B}$ that $\mathcal{A} = \mathcal{G}(a)$, $\mathcal{B} = \mathcal{H}(a)$, and $g^{(\mathcal{G})}(\mathcal{G}(a))h^{(\mathcal{H})}(\mathcal{H}(a)) \geq g^{(\mathcal{G})}(\mathcal{G}(1))h^{(\mathcal{H})}(\mathcal{H}(1))$. 
Since Lemma~\ref{complemma1} gives us $g^{(\mathcal{G})}(\mathcal{G}(a)) \leq g^{(\mathcal{G})}(\mathcal{G}(1))$ and $h^{(\mathcal{H})}(\mathcal{H}(a)) \leq h^{(\mathcal{H})}(\mathcal{H}(1))$, it follows that we actually have $g^{(\mathcal{G})}(\mathcal{G}(a))h^{(\mathcal{H})}(\mathcal{H}(a)) = g^{(\mathcal{G})}(\mathcal{G}(1))h^{(\mathcal{H})}(\mathcal{H}(1))$, $g^{(\mathcal{G})}(\mathcal{G}(a)) = g^{(\mathcal{G})}(\mathcal{G}(1))$ and $h^{(\mathcal{H})}(\mathcal{H}(a)) = h^{(\mathcal{H})}(\mathcal{H}(1))$.

Therefore, we will now assume that $\mathcal{A}$ and $\mathcal{B}$ are compressed.

Define $\mathcal{H}_0 = \{H \in \mathcal{H} \colon n \notin H\}$
and $\mathcal{H}_1 = \{H \backslash \{n\} \colon n \in H \in
\mathcal{H}\}$. Define $\mathcal{G}_0$, $\mathcal{G}_1$,
$\mathcal{A}_0$, $\mathcal{A}_1$, $\mathcal{B}_0$ and
$\mathcal{B}_1$ similarly. Since $\mathcal{A}$, $\mathcal{B}$,
$\mathcal{G}$ and $\mathcal{H}$ are compressed, we clearly have
that $\mathcal{A}_0$, $\mathcal{A}_1$, $\mathcal{B}_0$,
$\mathcal{B}_1$, $\mathcal{G}_0$, $\mathcal{G}_1$,
$\mathcal{H}_0$ and $\mathcal{H}_1$ are compressed. Since
$\mathcal{G}$ and $\mathcal{H}$ are hereditary, we clearly have
that $\mathcal{G}_0$, $\mathcal{G}_1$, $\mathcal{H}_0$ and
$\mathcal{H}_1$ are hereditary, $\mathcal{G}_1 \subseteq \mathcal{G}_0$ and $\mathcal{H}_1 \subseteq \mathcal{H}_0$. If $\mathcal{G}_1 = \emptyset$, then $\mathcal{G} \subseteq 2^{[m-1]}$ and hence we obtain the result immediately from the induction hypothesis. The same occurs if $\mathcal{H}_1 = \emptyset$. So we assume that $\mathcal{G}_1$ and $\mathcal{H}_1$ are non-empty. Since $\mathcal{G}_1 \subseteq \mathcal{G}_0$ and $\mathcal{H}_1 \subseteq \mathcal{H}_0$, $\mathcal{G}_0$ and $\mathcal{H}_0$ are non-empty too. Obviously, we have $\mathcal{A}_0 \subseteq \mathcal{G}_0 \subseteq 2^{[m-1]}$,
$\mathcal{A}_1 \subseteq \mathcal{G}_1 \subseteq 2^{[m-1]}$,
$\mathcal{B}_0 \subseteq \mathcal{H}_0 \subseteq 2^{[n-1]}$ and
$\mathcal{B}_1 \subseteq \mathcal{H}_1 \subseteq 2^{[n-1]}$.



Let $h_0 : \mathcal{H}_0 \rightarrow \mathbb{R}^+$ such that $h_0(H) = h(H)$ for each $H \in \mathcal{H}_0$. Let $h_1 : \mathcal{H}_1 \rightarrow \mathbb{R}^+$ such that $h_1(H) = h(H \cup \{n\})$ for each $H \in \mathcal{H}_1$ (note that $H \cup \{n\} \in \mathcal{H}(n)$ by definition of $\mathcal{H}_1$). By (a) and (b), we have the following consequences. For any $A, B \in \mathcal{H}_0$ with $A \subsetneq B$,
\begin{equation} h_0(A) = h(A) \geq 2h(B) = 2h_0(B). \label{main1}
\end{equation}
For any $C \in \mathcal{H}_0$ and any $i,j \in [n-1]$ with $i < j$,
\begin{equation} h_0(\delta_{i,j}(C)) = h(\delta_{i,j}(C)) \geq h(C) = h_0(C). \label{main2}
\end{equation}
For any $A, B \in \mathcal{H}_1$ with $A \subsetneq B$,
\begin{equation} h_1(A) = h(A \cup \{n\}) \geq 2h(B \cup \{n\}) = 2h_1(B). \label{main3}
\end{equation}
For any $C \in \mathcal{H}_1$ and any $i,j \in [n-1]$ with $i < j$,
\begin{equation} h_1(\delta_{i,j}(C)) = h(\delta_{i,j}(C) \cup \{n\}) = h(\delta_{i,j}(C \cup \{n\})) \geq h(C \cup \{n\}) = h_1(C). \label{main4}
\end{equation}
Thus, we have shown that properties (a) and (b) are inherited by $h_0$ and $h_1$.

Since $\mathcal{B} = \mathcal{B}_0 \cup \mathcal{B}(n)$, $\mathcal{B}_0 \cap \mathcal{B}(n) = \emptyset$ and $\mathcal{B}(n) = \{B \cup \{n\} \colon B \in \mathcal{B}_1\}$, we have
\begin{equation} h^{(\mathcal{H})}(\mathcal{B}) = {h}^{(\mathcal{H})}(\mathcal{B}_0) + h^{(\mathcal{H})}(\mathcal{B}(n)) = {h_0}^{(\mathcal{H}_0)}(\mathcal{B}_0) + {h_1}^{(\mathcal{H}_1)}(\mathcal{B}_1). \label{main4.2}
\end{equation}
Along the same lines,
\begin{align} h^{(\mathcal{H})}(\mathcal{H}(1)) &= h^{(\mathcal{H})}(\mathcal{H}_0(1)) + h^{(\mathcal{H})}(\{H \in \mathcal{H} \colon 1, n \in H\}) \nonumber \\
&= {h_0}^{(\mathcal{H}_0)}(\mathcal{H}_0(1)) + {h}^{(\mathcal{H})}(\{H \cup \{n\} \colon H \in \mathcal{H}_1(1)\}) \nonumber \\
&= {h_0}^{(\mathcal{H}_0)}(\mathcal{H}_0(1)) + {h_1}^{(\mathcal{H}_1)}(\mathcal{H}_1(1)). \label{main4.3}
\end{align}

Suppose $m < n$.  Clearly, $\mathcal{A}$ and $\mathcal{B}_0$ are
cross-intersecting. Since $m < n$, no set in $\mathcal{A}$ contains $n$, and hence $\mathcal{A}$ and $\mathcal{B}_1$ are cross-intersecting. Thus, by the induction hypothesis,
\begin{equation} g^{(\mathcal{G})}(\mathcal{A}){h_j}^{(\mathcal{H}_j)}(\mathcal{B}_j) \leq g^{(\mathcal{G})}(\mathcal{G}(1)){h_j}^{(\mathcal{H}_j)}(\mathcal{H}_j(1)) \quad \mbox{for each $j \in \{0, 1\}$}. \label{main4.4}
\end{equation}
Together with (\ref{main4.2}) and (\ref{main4.3}), this gives us
\begin{align} g^{(\mathcal{G})}(\mathcal{A}){h}^{(\mathcal{H})}(\mathcal{B})
&= g^{(\mathcal{G})}(\mathcal{A}){h_0}^{(\mathcal{H}_0)}(\mathcal{B}_0) + g^{(\mathcal{G})}(\mathcal{A}){h_1}^{(\mathcal{H}_1)}(\mathcal{B}_1) \nonumber \\
&\leq g^{(\mathcal{G})}(\mathcal{G}(1)){h_0}^{(\mathcal{H}_0)}(\mathcal{H}_0(1)) + g^{(\mathcal{G})}(\mathcal{G}(1)){h_1}^{(\mathcal{H}_1)}(\mathcal{H}_1(1)) \nonumber \\
&= g^{(\mathcal{G})}(\mathcal{G}(1)){h}^{(\mathcal{H})}(\mathcal{H}(1)). \nonumber 
\end{align}
By (\ref{main0.1}), equality holds throughout and hence $g^{(\mathcal{G})}(\mathcal{A}){h}^{(\mathcal{H})}(\mathcal{B}) = g^{(\mathcal{G})}(\mathcal{G}(1)){h}^{(\mathcal{H})}(\mathcal{H}(1))$. So in (\ref{main4.4}) we actually have equality. By the induction hypothesis, for each $j \in \{0,1\}$ we have $\mathcal{A} = \mathcal{G}(a_j)$ and $\mathcal{B}_j = \mathcal{H}_j(a_j)$ for some $a_j \in [m]$ such that $g^{(\mathcal{G})}(\mathcal{G}(a_j)) = g^{(\mathcal{G})}(\mathcal{G}(1))$ and ${h_j}^{(\mathcal{H}_j)}(\mathcal{H}_j(a_j)) = {h_j}^{(\mathcal{H}_j)}(\mathcal{H}_j(1))$. So $g^{(\mathcal{G})}(\mathcal{G}(a_0)) > 0$ and hence $\mathcal{G}(a_0) \neq \emptyset$. Thus, since $\mathcal{G}$ is hereditary, $\{a_0\} \in \mathcal{A}$. Since $\mathcal{A}$ and $\mathcal{B}$ are cross-intersecting, $\mathcal{B} \subseteq \mathcal{H}(a_0)$. Since $\mathcal{G}(a_0)$ and $\mathcal{H}(a_0)$ are cross-intersecting, it follows by the choice of $\mathcal{A}$ and $\mathcal{B}$ that $\mathcal{A} = \mathcal{G}(a_0)$ and $\mathcal{B} = \mathcal{H}(a_0)$. Thus, since $g^{(\mathcal{G})}(\mathcal{A}){h}^{(\mathcal{H})}(\mathcal{B}) = g^{(\mathcal{G})}(\mathcal{G}(1)){h}^{(\mathcal{H})}(\mathcal{H}(1))$, and since Lemma~\ref{complemma1} gives us $g^{(\mathcal{G})}(\mathcal{G}(a_0)) \leq g^{(\mathcal{G})}(\mathcal{G}(1))$ and $h^{(\mathcal{H})}(\mathcal{H}(a_0)) \leq h^{(\mathcal{H})}(\mathcal{H}(1))$, we have $g^{(\mathcal{G})}(\mathcal{G}(a_0)) = g^{(\mathcal{G})}(\mathcal{G}(1))$ and $h^{(\mathcal{H})}(\mathcal{H}(a_0)) = h^{(\mathcal{H})}(\mathcal{H}(1))$.\medskip

Now suppose $m=n$. Similarly to $h_0$ and $h_1$, let $g_0 : \mathcal{G}_0 \rightarrow \mathbb{R}^+$ such that $g_0(G) = g(G)$ for each $G \in \mathcal{G}_0$, and let $g_1 : \mathcal{G}_1 \rightarrow \mathbb{R}^+$ such that $g_1(G) = g(G \cup \{n\})$ for each $G \in \mathcal{G}_1$ (note that since $m = n$, $G \cup \{n\} \in \mathcal{G}(n)$ by definition of $\mathcal{G}_1$). Then the properties (a) and (b) are inherited by $g_0$ and $g_1$ in the same way they are inherited by $h_0$ and $h_1$ as shown above; that is, similarly to (\ref{main1})--(\ref{main4}), we have the following. For any $A, B \in \mathcal{G}_0$ with $A \subsetneq B$,
\begin{equation} g_0(A) \geq 2g_0(B). \label{main5}
\end{equation}
For any $C \in \mathcal{G}_0$ and any $i,j \in [n-1]$ with $i < j$,
\begin{equation} g_0(\delta_{i,j}(C)) \geq g_0(C). \label{main6}
\end{equation}
For any $A, B \in \mathcal{G}_1$ with $A \subsetneq B$,
\begin{equation} g_1(A) \geq 2g_1(B). \label{main7}
\end{equation}
For any $C \in \mathcal{G}_1$ and any $i,j \in [n-1]$ with $i < j$,
\begin{equation} g_1(\delta_{i,j}(C)) \geq g_1(C). \label{main8}
\end{equation}

Similarly to (\ref{main4.2}) and (\ref{main4.3}), we have
\begin{gather} g^{(\mathcal{G})}(\mathcal{A}) = {g_0}^{(\mathcal{G}_0)}(\mathcal{A}_0) + {g_1}^{(\mathcal{G}_1)}(\mathcal{A}_1), \label{main4.1} \\
g^{(\mathcal{G})}(\mathcal{G}(1)) = {g_0}^{(\mathcal{G}_0)}(\mathcal{G}_0(1)) + {g_1}^{(\mathcal{G}_1)}(\mathcal{G}_1(1)). \label{main4.6}
\end{gather}

Clearly, $\mathcal{A}_0$ and $\mathcal{B}_0$ are cross-intersecting, and since $n = m$, so are $\mathcal{A}_0$ and $\mathcal{B}_1$, and also $\mathcal{A}_1$ and $\mathcal{B}_0$.

Let us first assume that $\mathcal{A}_1$ and $\mathcal{B}_1$ are cross-intersecting too. Then, by the induction hypothesis,
\begin{equation} {g_i}^{(\mathcal{G}_i)}(\mathcal{A}_i){h_j}^{(\mathcal{H}_j)}(\mathcal{B}_j) \leq {g_i}^{(\mathcal{G}_i)}(\mathcal{G}_i(1)){h_j}^{(\mathcal{H}_j)}(\mathcal{H}_j(1)) \quad \mbox{for any $i, j \in \{0, 1\}$.} \label{main4.5}
\end{equation}
Together with (\ref{main4.2}), (\ref{main4.3}), (\ref{main4.1}) and (\ref{main4.6}), this gives us
\begin{align} g^{(\mathcal{G})}(\mathcal{A}){h}^{(\mathcal{H})}(\mathcal{B}) &= {g_0}^{(\mathcal{G}_0)}(\mathcal{A}_0) {h_0}^{(\mathcal{H}_0)}(\mathcal{B}_0) + {g_0}^{(\mathcal{G}_0)}(\mathcal{A}_0) {h_1}^{(\mathcal{H}_1)}(\mathcal{B}_1) + \nonumber \\
& \quad \mbox{ } {g_1}^{(\mathcal{G}_1)}(\mathcal{A}_1) {h_0}^{(\mathcal{H}_0)}(\mathcal{B}_0) + {g_1}^{(\mathcal{G}_1)}(\mathcal{A}_1){h_1}^{(\mathcal{H}_1)}(\mathcal{B}_1) \nonumber \\
&\leq {g_0}^{(\mathcal{G}_0)}(\mathcal{G}_0(1)){h_0}^{(\mathcal{H}_0)}(\mathcal{H}_0(1)) + {g_0}^{(\mathcal{G}_0)}(\mathcal{G}_0(1)){h_1}^{(\mathcal{H}_1)}(\mathcal{H}_1(1)) + \nonumber \\
& \quad \mbox{ } {g_1}^{(\mathcal{G}_1)}(\mathcal{G}_1(1)){h_0}^{(\mathcal{H}_0)}(\mathcal{H}_0(1)) + {g_1}^{(\mathcal{G}_1)}(\mathcal{G}_1(1)){h_1}^{(\mathcal{H}_1)}(\mathcal{H}_1(1)) \nonumber \\
&= g^{(\mathcal{G})}(\mathcal{G}(1)){h}^{(\mathcal{H})}(\mathcal{H}(1)). \nonumber
\end{align}
By (\ref{main0.1}), equality holds throughout and hence $g^{(\mathcal{G})}(\mathcal{A}){h}^{(\mathcal{H})}(\mathcal{B}) = g^{(\mathcal{G})}(\mathcal{G}(1)){h}^{(\mathcal{H})}(\mathcal{H}(1))$. So in (\ref{main4.5}) we actually have equality. By the induction hypothesis, we particularly have $\mathcal{A}_0 = \mathcal{G}_0(a_0)$ and $\mathcal{B}_0 = \mathcal{H}_0(a_0)$ for some $a_0 \in [n-1]$ such that ${g_0}^{(\mathcal{G}_0)}(\mathcal{G}_0(a_0)) = {g_0}^{(\mathcal{G}_0)}(\mathcal{G}_0(1))$ and ${h_0}^{(\mathcal{H}_0)}(\mathcal{H}_0(a_0)) = {h_0}^{(\mathcal{H}_0)}(\mathcal{H}_0(1))$. Recall that $\{1\} \in \mathcal{G}$. So $\{1\} \in \mathcal{G}_0$ and hence ${g_0}^{(\mathcal{G}_0)}(\mathcal{G}_0(1)) > 0$. So ${g_0}^{(\mathcal{G}_0)}(\mathcal{G}_0(a_0)) > 0$ and hence $\mathcal{G}_0(a_0) \neq \emptyset$. Thus, since $\mathcal{G}$ is hereditary, $\{a_0\} \in \mathcal{A}$. Since $\mathcal{A}_0$ and $\mathcal{B}$ are cross-intersecting, $\mathcal{B} \subseteq \mathcal{H}(a_0)$. Similarly, we obtain $\mathcal{A} \subseteq \mathcal{G}(a_0)$. As in the case $m < n$, we conclude that $\mathcal{A} = \mathcal{G}(a_0)$, $\mathcal{B} = \mathcal{H}(a_0)$, $g^{(\mathcal{G})}(\mathcal{G}(a_0)) = g^{(\mathcal{G})}(\mathcal{G}(1))$ and $h^{(\mathcal{H})}(\mathcal{H}(a_0)) = h^{(\mathcal{H})}(\mathcal{H}(1))$.\medskip

We will now show that indeed $\mathcal{A}_1$ and $\mathcal{B}_1$ are cross-intersecting. Suppose they are not. Then there exists $A_1 \in \mathcal{A}_1$ such that $A_1 \cap B = \emptyset$ for some $B \in \mathcal{B}_1$. Let $B_1 = [n-1] \backslash A_1$, $A_1' = A_1 \cup \{n\}$, $B_1' = B_1 \cup \{n\}$. Since $A_1 \in \mathcal{A}_1$, $A_1' \in \mathcal{A}$.

If $A_1 = [n-1]$, then $B = B_1$. Suppose that $A_1 \neq [n-1]$ and $B \neq B_1$. Then $B \subsetneq [n-1] \backslash A_1$ and hence $[n-1] \backslash (A_1 \cup B) \neq \emptyset$. Let $c \in [n-1] \backslash (A_1 \cup B)$. Since $B \in \mathcal{B}_1$, $B \cup \{n\} \in \mathcal{B}$. Let $C = \delta_{c,n}(B \cup \{n\})$. Since $c \notin B \cup \{n\}$, $C = B \cup \{c\}$. Since $\mathcal{B}$ is compressed, $C \in \mathcal{B}$. However, since $c \notin A_1'$ and $A_1 \cap B = \emptyset$, we have $A_1' \cap C = \emptyset$, which is a contradiction as $\mathcal{A}$ and $\mathcal{B}$ are cross-intersecting. 

We have therefore shown that
\begin{equation}\mbox{$B_1$ is the unique set in $\mathcal{B}_1$ that does not intersect $A_1$.} \label{main8.1}
\end{equation}
%
By a similar argument,
\begin{equation}\mbox{$A_1$ is the unique set in $\mathcal{A}_1$ that does not intersect $B_1$.} \label{main8.2}
\end{equation}
Since $B_1 \in \mathcal{B}_1$, $B_1' \in \mathcal{B}$. Since $\mathcal{A}$ and $\mathcal{B}$ are compressed, 
\begin{equation} \mbox{ for any $p \in [n-1]$, $\delta_{p,n}(A_1') \in \mathcal{A}$ and $\delta_{p,n}(B_1') \in \mathcal{B}$.} \label{main9}
\end{equation}

Since $A_1 \cap B_1' = A_1 \cap B_1 = \emptyset$ and $B_1 \cap A_1' = B_1 \cap A_1 = \emptyset$, $A_1 \notin \mathcal{A}$ and $B_1 \notin \mathcal{B}$. Let $\mathcal{A}' = \mathcal{A} \cup \{A_1\}$, $\mathcal{A}'' = \mathcal{A} \backslash \{A_1'\}$, $\mathcal{B}' = \mathcal{B} \backslash \{B_1'\}$, $\mathcal{B}'' = \mathcal{B} \cup \{B_1\}$. By (\ref{main8.1}), $\mathcal{A}'$ and $\mathcal{B}'$ are cross-intersecting. By (\ref{main8.2}), $\mathcal{A}''$ and $\mathcal{B}''$ are cross-intersecting. Since $\mathcal{G}$ and $\mathcal{H}$ are hereditary, and since $A_1' \in \mathcal{A} \subseteq \mathcal{G}$ and $B_1' \in \mathcal{B} \subseteq \mathcal{H}$, we have $A_1 \in \mathcal{G}$ and $B_1 \in \mathcal{H}$, and hence $\mathcal{A}', \mathcal{A}'' \in \mathcal{G}$ and $\mathcal{B}', \mathcal{B}'' \in \mathcal{H}$. 

Let $x = g^{(\mathcal{G})}(\mathcal{A})$ and $x_1 = g(A_1')$. Let $y = h^{(\mathcal{H})}(\mathcal{B})$ and $y_1 = h(B_1')$. We have
\begin{align} &g^{(\mathcal{G})}(\mathcal{A}') = x + g(A_1) \geq x + 2g(A_1') = x + 2x_1,  \nonumber \\
&g^{(\mathcal{G})}(\mathcal{A}'') = x - g(A_1') = x - x_1, \nonumber \\
&h^{(\mathcal{H})}(\mathcal{B}') = y - h(B_1') = y - y_1, \nonumber \\
&{h}^{(\mathcal{H})}(\mathcal{B}'') = y + h(B_1) \geq y + 2h(B_1') = y + 2y_1.  \nonumber
\end{align}
By the choice of $\mathcal{A}$ and $\mathcal{B}$, 
\[g^{(\mathcal{G})}(\mathcal{A}'){h}^{(\mathcal{H})}(\mathcal{B}') \leq g^{(\mathcal{G})}(\mathcal{A}){h}^{(\mathcal{H})}(\mathcal{B}) \quad \mbox{and} \quad g^{(\mathcal{G})}(\mathcal{A}''){h}^{(\mathcal{H})}(\mathcal{B}'') \leq g^{(\mathcal{G})}(\mathcal{A}){h}^{(\mathcal{H})}(\mathcal{B}).\]
So we have
\begin{align} &(x + 2x_1)(y - y_1) \leq xy \quad \mbox{and} \quad (x - x_1)(y + 2y_1) \leq xy \nonumber \\
&\Rightarrow 2x_1y \leq xy_1 + 2x_1y_1 \quad \mbox{and} \quad 2y_1x \leq x_1y + 2x_1y_1 \nonumber \\
&\Rightarrow  2x_1y + 2y_1x \leq (y_1x + 2x_1y_1) + (x_1y + 2x_1y_1) \nonumber \\
&\Rightarrow  x_1y + y_1x \leq 4x_1y_1. \nonumber
\end{align}

Suppose that $A_1 = \emptyset$. Then $A_1' = \{n\}$ and $B_1' = [n]$. By (\ref{main9}), the sets $\{1\}, \dots, \{n\}$ are all in $\mathcal{A}$, and obviously no proper subset of $[n]$ intersects each of these sets. Thus, by the cross-intersection condition, $B_1'$ is the only set that is in $\mathcal{B}$. So ${h}^{(\mathcal{H})}(\mathcal{B}'') = h(B_1') + h(B_1) \geq h(B_1') + 2h(B_1') = 3h(B_1') = 3{h}^{(\mathcal{H})}(\mathcal{B})$. Since $\{1\}, \{n\} \in \mathcal{A}$ and $A_1' = \{n\}$, we have $2g(A_1') \leq g(A_1') + g(\delta_{1,n}(A_1')) = g(\{n\}) + g(\{1\}) \leq g^{(\mathcal{G})}(\mathcal{A})$ and hence $g^{(\mathcal{G})}(\mathcal{A}'') = g^{(\mathcal{G})}(\mathcal{A}) - g(A_1') \geq g^{(\mathcal{G})}(\mathcal{A})/2$. So we obtain $g^{(\mathcal{G})}(\mathcal{A}'') {h}^{(\mathcal{H})}(\mathcal{B}'') \geq (g^{(\mathcal{G})}(\mathcal{A})/2) (3{h}^{(\mathcal{H})}(\mathcal{B}))$, which contradicts $g^{(\mathcal{G})}(\mathcal{A}''){h}^{(\mathcal{H})}(\mathcal{B}'') \leq g^{(\mathcal{G})}(\mathcal{A}) {h}^{(\mathcal{H})}(\mathcal{B})$. 

Therefore, $A_1 \neq \emptyset$. Also, $B_1 \neq \emptyset$ because otherwise, by the same argument, we obtain $g^{(\mathcal{G})}(\mathcal{A}'){h}^{(\mathcal{H})}(\mathcal{B}') \geq (3g^{(\mathcal{G})}(\mathcal{A}))({h}^{(\mathcal{H})}(\mathcal{B})/2)$, which contradicts $g^{(\mathcal{G})}(\mathcal{A}'){h}^{(\mathcal{H})}(\mathcal{B}') \leq g^{(\mathcal{G})}(\mathcal{A}){h}^{(\mathcal{H})}(\mathcal{B})$. Since $A_1 \neq \emptyset \neq B_1$, it follows by definition of $B_1$ that  $[n-1] \backslash A_1 \neq \emptyset$ and $[n-1] \backslash B_1 \neq \emptyset$. Let $a \in [n-1] \backslash A_1$ and $b \in [n-1] \backslash B_1$. Let $A_1'' = \delta_{a,n}(A_1')$ and $B_1'' = \delta_{b,n}(B_1')$. So $A_1'' \neq A_1'$ and $B_1'' \neq B_1'$. By (\ref{main9}), $A_1'' \in \mathcal{A}$ and $B_1'' \in \mathcal{B}$. By (b), $g(A_1'') \geq g(A_1')$ and $h(B_1'') \geq h(B_1')$. We therefore have $x \geq x_1 + g(A_1'') \geq 2x_1$ and $y \geq y_1 + h(B_1'') \geq 2y_1$. So $x_1y + y_1x \geq x_1(2y_1) + y_1(2x_1) = 4x_1y_1$. Together with $x_1y + y_1x \leq 4x_1y_1$, this gives us $x_1y + y_1x = 4x_1y_1$. Consequently, we have $\mathcal{A} = \{A_1', A_1''\}$, $\mathcal{B} = \{B_1', B_1''\}$, $g(A_1'') = g(A_1')$ and $h(B_1'') = h(B_1')$. Let $A_2 = [|A_1'|]$, $B_2 = [|B_1'|]$ and $I = \{1\}$. Since $\mathcal{G}$ is compressed and $A_1' \in \mathcal{A} \subseteq \mathcal{G}$, $A_2 \in \mathcal{G}$ and $g(A_2) \geq g(A_1')$. Similarly, $B_2 \in \mathcal{H}$ and $h(B_2) \geq h(B_1')$. Since $A_1 \neq \emptyset$, we have $|A_1'| \geq 2$ and hence $|A_2| \geq 2$. Since $\mathcal{G}$ is hereditary and $I \subset A_2 \in \mathcal{G}$, $I \in \mathcal{G}$ and $g(I) \geq 2g(A_2)$. Similarly, $I \in \mathcal{H}$ and $h(I) \geq 2h(B_2)$. Let $\mathcal{C} = \{I,A_2\}$ and $\mathcal{D} = \{I,B_2\}$. So $\mathcal{C} \subset \mathcal{G}$ and $\mathcal{D} \subset \mathcal{H}$. Also, $\mathcal{C}$ and $\mathcal{D}$ are cross-intersecting. We have 
\begin{align} g^{(\mathcal{G})}(\mathcal{C}) h^{(\mathcal{H})}(\mathcal{D}) &= (g(I) + g(A_2))(h(I) + h(B_2)) \geq (3g(A_2))(3h(B_2)) 
= 9g(A_2)h(B_2) \nonumber \\
&\geq 9g(A_1')h(B_1') = 9\left(\frac{g(A_1') + g(A_1'')}{2}\right)\left(\frac{h(B_1') + h(B_1'')}{2}\right) \nonumber \\
&= \frac{9}{4}g^{(\mathcal{G})}(\mathcal{A}) h^{(\mathcal{H})}(\mathcal{B}), \nonumber
\end{align}
which contradicts the choice of $\mathcal{A}$ and $\mathcal{B}$.~\hfill{$\Box$}

\section{Proofs of Theorems~\ref{main} and \ref{mainpartial}} \label{Proofmain}

We now define a compression operation for labeled sets. For any $x, y \in \mathbb{N}$, let
\[ \gamma_{x,y}(A) = \left\{ \begin{array}{ll}
(A \backslash \{(x,y)\}) \cup \{(x,1)\} & \mbox{if $(x,y) \in A$};\\
A & \mbox{otherwise}
\end{array} \right. \]
for any labeled set $A$, and let
\[\Gamma_{x,y}(\mathcal{A}) = \{\gamma_{x,y}(A) \colon A \in \mathcal{A}, \gamma_{x,y}(A) \notin \mathcal{A}\} \cup \{A \in \mathcal{A} \colon \gamma_{x,y}(A) \in
\mathcal{A}\}\]
for any family $\mathcal{A}$ of labeled sets.

Note that $|\Gamma_{x,y}(\mathcal{A})| = |\mathcal{A}|$ and that if $\mathcal{A} \subseteq \mathcal{L}_{\bf c}^{(r)}$, then $\Gamma_{x,y}(\mathcal{A}) \subseteq \mathcal{L}_{\bf c}^{(r)}$. It is easy to check that if $\mathcal{A}$ and $\mathcal{B}$ are cross-intersecting families of labeled sets, then so are $\Gamma_{x,y}(\mathcal{A})$ and $\Gamma_{x,y}(\mathcal{B})$. We prove more than this.

The \emph{Cartesian product} of two sets $X$ and $Y$, denoted by $X \times Y$, is the set $\{(x, y) \colon x \in X, \, y \in Y\}$. For any IP sequence ${\bf c} = (c_1, \dots, c_n)$ and any $r \in [n]$, let $\mathcal{L}_{\bf c}^{(\leq r)}$ denote the union $\bigcup_{i = 1}^r \mathcal{L}_{\bf c}^{(i)}$.

\begin{lemma}\label{gamma} Let ${\bf c} = (c_1, \dots, c_m)$ and ${\bf d} = (d_1, \dots, d_n)$ be IP sequences. Let $x,y \in \mathbb{N}$, $y \geq 2$. Let $l = \max\{m,n\}$ and $h = \max\{c_m, d_n\}$. Let $V \subseteq [l] \times [2,h]$. Let $\mathcal{A} \subseteq \mathcal{L}_{\bf c}^{(\leq m)}$ and $\mathcal{B} \subseteq \mathcal{L}_{\bf d}^{(\leq n)}$ such that $(A \cap B) \backslash V \neq \emptyset$ for every $A \in \mathcal{A}$ and every $B \in \mathcal{B}$. Then $(C \cap D) \backslash (V \cup \{(x,y)\}) \neq \emptyset$ for every $C \in \Gamma_{x,y}(\mathcal{A})$ and every $D \in \Gamma_{x,y}(\mathcal{B})$.
\end{lemma}
\textbf{Proof.} Let $C \in \Gamma_{x,y}(\mathcal{A})$ and $D \in \Gamma_{x,y}(\mathcal{B})$. We first show that $(C \cap D) \backslash
V \neq \emptyset$. Let $C' = (C \backslash \{(x,1)\}) \cup \{(x,y)\}$ and $D' = (D \backslash \{(x,1)\}) \cup \{(x,y)\}$. If $C \in \mathcal{A}$ and $D \in \mathcal{B}$, then $(C \cap D) \backslash V \neq \emptyset$. If $C \notin \mathcal{A}$ and $D \notin \mathcal{B}$, then $(x,1)$ is in both $C$ and $D$, and hence, since $(x,1) \notin V$, $(x,1) \in (C \cap D) \backslash V$. Suppose $C \notin \mathcal{A}$ and $D \in \mathcal{B}$.  So $(x,1) \in C$ and $C' \in \mathcal{A}$. If $(x,y) \notin D$, then, since $C' \in \mathcal{A}$ and $D \in \mathcal{B}$, $0 < |(C' \cap D) \backslash V| \leq |(C \cap D) \backslash V|$. If $(x,y) \in D$, then $\gamma_{x,y}(D) \in \mathcal{B}$ (because otherwise $D \notin \Gamma_{x,y}(\mathcal{B})$), and hence, since $C' \in \mathcal{A}$, $0 < |(C' \cap \gamma_{x,y}(D)) \backslash V| = |(C \cap D) \backslash V|$. Similarly, if $C \in \mathcal{A}$ and $D \notin \mathcal{B}$, then $(C \cap D) \backslash
V \neq \emptyset$.

Now suppose $(C \cap D) \backslash (V \cup \{(x,y)\}) = \emptyset$.
Since $(C \cap D) \backslash V \neq \emptyset$, $(x,y) \in C \cap D$. So $C, \gamma_{x,y}(C) \in \mathcal{A}$, $D, \gamma_{x,y}(D) \in \mathcal{B}$ and $|(C \cap \gamma_{x,y}(D)) \backslash V| = |(C \cap D) \backslash (V \cup \{(x,y)\})| = 0$, a contradiction.~\hfill{$\Box$}

\begin{cor}\label{deltacor} Let ${\bf c} = (c_1, \dots, c_m), {\bf d} = (d_1, \dots, d_n), h$ and $l$ be as in Lemma~\ref{gamma}. Let $\mathcal{A} \subseteq \mathcal{L}_{\bf c}^{(\leq m)}$ and $\mathcal{B} \subseteq \mathcal{L}_{\bf d}^{(\leq n)}$ such that $\mathcal{A}$ and $\mathcal{B}$ are cross-intersecting. Let
\begin{align} \mathcal{A}^* &= \Gamma_{l,h} \circ \dots \circ \Gamma_{l,2} \circ \dots \circ \Gamma_{2,h} \circ \dots \circ \Gamma_{2,2} \circ \Gamma_{1,h} \circ \dots \circ \Gamma_{1,2}(\mathcal{A}),\nonumber \\
\mathcal{B}^* &= \Gamma_{l,h} \circ \dots \circ \Gamma_{l,2} \circ \dots \circ \Gamma_{2,h} \circ \dots \circ \Gamma_{2,2} \circ \Gamma_{1,h} \circ \dots \circ \Gamma_{1,2}(\mathcal{B}). \nonumber
\end{align}
Then $A \cap B \cap ([l] \times [1]) \neq \emptyset$ for any $A \in \mathcal{A}^*$ and any $B \in \mathcal{B}^*$.
\end{cor}
\textbf{Proof.} Let $Z = [l] \times [2,h]$. By repeated application of Lemma~\ref{gamma},
$(A \cap B) \backslash Z \neq \emptyset$ for any $A \in \mathcal{A}^*$ and any $B \in \mathcal{B}^*$. The result follows since $(A \cap B) \backslash Z = A \cap B \cap ([l] \times [1])$.~\hfill{$\Box$}\\

The next lemma is needed for the characterisation of the extremal structures in Theorems~\ref{main} and \ref{mainpartial}. Note that according to the notation in Section~\ref{Weightedsection}, $\mathcal{L}_{\bf c}^{(r)}((a,b)) = \{A \in \mathcal{L}_{\bf c}^{(r)} \colon (a,b) \in A\}$.

\begin{lemma}\label{ss_cross_int} Let ${\bf c} = (c_1, \dots, c_m), {\bf d} = (d_1, \dots, d_n), h$ and $l$ be as in Lemma~\ref{gamma}. Suppose that $c_1 \geq 3$ and $d_1 \geq 3$. Let $r \in [m]$ and $s \in [n]$. Let $\mathcal{A} \subseteq \mathcal{L}_{\bf c}^{(r)}$ and $\mathcal{B} \subseteq \mathcal{L}_{\bf d}^{(s)}$ such that $\mathcal{A}$ and $\mathcal{B}$ are cross-intersecting. Suppose that, for some $(x,y), (u,v) \in [l] \times [h]$, we have $\Gamma_{x,y}(\mathcal{A}) = \mathcal{L}_{\bf c}^{(r)}((u,v))$, $\Gamma_{x,y}(\mathcal{B}) = \mathcal{L}_{\bf d}^{(s)}((u,v))$, and either $\mathcal{A} \neq
\Gamma_{x,y}(\mathcal{A})$ or $\mathcal{B} \neq \Gamma_{x,y}(\mathcal{B})$. Then $\mathcal{A} = \mathcal{L}_{\bf c}^{(r)}((x,y))$ and $\mathcal{B} = \mathcal{L}_{\bf d}^{(s)}((x,y))$.
\end{lemma}
\textbf{Proof.} We may assume that $\mathcal{A} \neq \Gamma_{x,y}(\mathcal{A})$. So there exists $A_1 \in \mathcal{A} \backslash \Gamma_{x,y}(\mathcal{A})$ such that $\gamma_{x,y}(A_1) \in \Gamma_{x,y}(\mathcal{A}) \backslash \mathcal{A}$. Let $A_1' = \gamma_{x,y}(A_1)$. Thus, $(x,y) \in A_1$ and $(u,v) \in A_1' = (A_1 \backslash \{(x,y)\}) \cup \{(x,1)\}$. 

Suppose that $(u,v) \neq (x,1)$. Then $(u,v) \in A_1$. So $A_1 \in \mathcal{L}_{\bf c}^{(r)}((u,v))$ and hence $A_1 \in \Gamma_{x,y}(\mathcal{A})$, a contradiction. 

Therefore, $(u,v) = (x,1)$. Since $A_1 \neq A_1'$, $(x,y) \neq (x,1)$.

Let $A^* \in \mathcal{L}_{\bf c}^{(r)}((x,y))$. Let $x_1, \dots, x_{s-1}$ be distinct elements of $[n] \backslash \{x\}$. We are given that $3 \leq d_1 \leq \dots d_n$. By definition of a labeled set, for each $i \in [n]$ we have $|A \cap [d_i]| \leq 1$ for all $A \in \mathcal{L}_{\bf c}^{(r)}$. Thus, $|[d_i] \backslash (A_1 \cup A^*)| \geq d_i - 2 \geq 1$ for each $i \in [n]$. For each $i \in [s-1]$, let $y_i \in [d_i] \backslash (A_1 \cup A^*)$. Let $B^* = \{(x,y), (x_1,y_1), \dots, (x_{s-1},y_{s-1})\}$. So $B^* \in \mathcal{L}_{\bf d}^{(s)}((x,y))$. Since $\Gamma_{x,y}(\mathcal{B}) = \mathcal{L}_{\bf d}^{(s)}((x,1))$, either $B^* \in \mathcal{B}$ or $\gamma_{x,y}(B^*) \in \mathcal{B}$. However, $\gamma_{x,y}(B^*) \cap A_1 = \emptyset$. So $B^* \in \mathcal{B}$. Since $\Gamma_{x,y}(\mathcal{A}) = \mathcal{L}_{\bf c}^{(r)}((x,1))$, either $A^* \in \mathcal{A}$ or $\gamma_{x,y}(A^*) \in \mathcal{A}$. However, $\gamma_{x,y}(A^*) \cap B^* = \emptyset$. So $A^* \in \mathcal{A}$. 

We have therefore shown that $\mathcal{L}_{\bf c}^{(r)}((x,y)) \subseteq \mathcal{A}$. Since $|\Gamma_{x,y}(\mathcal{A})| = |\mathcal{L}_{\bf c}^{(r)}((x,1))| = |\mathcal{L}_{\bf c}^{(r)}((x,y))|$, we actually have $\mathcal{A} = \mathcal{L}_{\bf c}^{(r)}((x,y))$. Clearly, for each $L \in \mathcal{L}_{\bf d}^{(s)}$ with $(x,y) \notin L$, there exists $L' \in \mathcal{L}_{\bf c}^{(r)}((x,y))$ such that $L \cap L' = \emptyset$. Thus, since $\mathcal{A} = \mathcal{L}_{\bf c}^{(r)}((x,y))$, each set in $\mathcal{B}$ must contain $(x,y)$. So $\mathcal{B} \subseteq \mathcal{L}_{\bf d}^{(s)}((x,y))$. Since $|\Gamma_{x,y}(\mathcal{B})| = |\mathcal{L}_{\bf d}^{(s)}((x,1))| = |\mathcal{L}_{\bf d}^{(s)}((x,y))|$, we actually have $\mathcal{B} = \mathcal{L}_{\bf d}^{(s)}((x,y))$.~\hfill{$\Box$}\\

Let $\delta_{i,j}$ be as in Section~\ref{Compsection}. For any set $X$ and any $r \in \mathbb{N}$, let ${X \choose \leq r}$ denote the family $\{A \subseteq X \colon |A| \leq r\}$.

\begin{lemma}\label{lemmaweighted} 
Let ${\bf c}$ be an IP sequence $(c_1, \dots, c_n)$ such that $c_1 \geq 3$. Let $r \in [n]$. Let $w \colon {[n] \choose \leq r} \rightarrow \mathbb{N}$ such that for each $A \in {[n] \choose \leq r}$,
\[w(A) = \left| \left\{ L \in \mathcal{L}_{\bf c}^{(r)} \colon L \cap ([n] \times [1]) = A \times [1] \right\} \right|.\]
Then: \\
(i) $w(A) \geq 2w(A')$ for any $A, A' \in {[n] \choose \leq r}$ with $A \subsetneq A'$. \\
(ii) $w(\delta_{i,j}(A)) \geq w(A)$ for any $A \in {[n] \choose \leq r}$ and any $i,j \in [n]$ with $i < j$.
\end{lemma}
\textbf{Proof of Lemma \ref{lemmaweighted}.} 
(i) Let $A, A' \in {[n] \choose \leq r}$ with $A \subsetneq A'$. Let $B = A' \backslash A$. So $|B| \geq 1$. For each $L \in \mathcal{L}_{\bf c}^{(r)}$, let $\sigma(L) = \{i \in [n] \colon (i,a) \in L \mbox{ for some } a \in [c_i]\}$. We have
\begin{align} w(A) &\geq \left| \left\{ L \in \mathcal{L}_{\bf c}^{(r)} \colon L \cap ([n] \times [1]) = A \times [1], \, B \subset \sigma(L) \right\} \right| \nonumber \\
&= \sum_{E \in {[n] \backslash (A \cup B) \choose r - |A| - |B|}} \prod_{b \in B} (c_b - 1) \prod_{e \in E} (c_e - 1) = \prod_{b \in B} (c_b - 1) \left( \sum_{E \in {[n] \backslash A' \choose r - |A'|}} \prod_{e \in E} (c_e - 1) \right) \nonumber \\
&= w(A') \prod_{b \in B} (c_b - 1) \geq 2^{|B|}w(A') \geq 2w(A'). \nonumber
\end{align}
(ii) Let $A \in {[n] \choose \leq r}$, and let $i,j \in [n]$ with $i < j$. Suppose $\delta_{i,j}(A) \neq A$. Then $j \in A$, $i \notin A$ and $\delta_{i,j}(A) = (A \backslash \{j\}) \cup \{i\}$. Let $B = A \backslash \{j\}$, $\mathcal{E}_0 = {[n] \backslash (B \cup \{i, j\}) \choose r - |A|}$, $\mathcal{E}_1 = \left\{E \in {[n] \backslash (B \cup \{i\}) \choose r - |A|} \colon j \in E \right\}$, $\mathcal{E}_2 = \left \{E \in {[n] \backslash (B \cup \{j\}) \choose r - |A|} \colon i \in E \right\}$.  We have
\begin{align} w(B \cup \{i\}) &= \sum_{E \in {[n] \backslash (B \cup \{i\}) \choose r - |A|}}\prod_{e \in E}(c_e - 1) = \sum_{D \in \mathcal{E}_0}\prod_{d \in D}(c_d - 1) + \sum_{F \in \mathcal{E}_1}\prod_{f \in F}(c_f - 1) \nonumber \\
&\geq \sum_{D \in \mathcal{E}_0}\prod_{d \in D}(c_d - 1) + \sum_{F \in \mathcal{E}_1}\prod_{f \in F}(c_f - 1)\frac{c_i-1}{c_j-1} \quad \quad \mbox{(since $c_i \leq c_j$)} \nonumber \\
&= \sum_{D \in \mathcal{E}_0}\prod_{d \in D}(c_d - 1) + \sum_{F \in \mathcal{E}_2}\prod_{f \in F}(c_f - 1) = \sum_{E \in {[n] \backslash (B \cup \{j\}) \choose r - |A|}}\prod_{e \in E}(c_e - 1) \nonumber \\
&= w(B \cup \{j\}) \nonumber
\end{align}
and hence $w(\delta_{i,j}(A)) \geq w(A)$.~\hfill{$\Box$}
\\

We now prove Theorem~\ref{mainpartial}, and then we prove Theorem~\ref{main}.\\
\\
\textbf{Proof of Theorem~\ref{mainpartial}.} Let $\mathcal{X} = \{X \in \mathcal{L}_{\bf c}^{(r)} \colon (1,1) \in X\}$ and $\mathcal{Y} = \{Y \in \mathcal{L}_{\bf c}^{(s)} \colon (1,1) \in Y\}$. Let $f : S_{\bf c}^{(r)} \rightarrow \mathcal{L}_{\bf c}^{(r)}$ such that for each ${\bf a} = (a_1, \dots, a_m)$ in $S_{\bf c}^{(r)}$, $f({\bf a}) = \{(i,a_i) \colon i \in [m], \, a_i \neq 0\}$. Let $g : S_{\bf d}^{(s)} \rightarrow \mathcal{L}_{\bf d}^{(s)}$ such that for each ${\bf b} = (b_1, \dots, b_n)$ in $S_{\bf d}^{(s)}$, $g({\bf b}) = \{(j,b_j) \colon j \in [n], \, b_j \neq 0\}$. Let $\mathcal{A} = \{f({\bf a}) \colon {\bf a} \in A\}$ and $\mathcal{B} = \{g({\bf b}) \colon {\bf b} \in B\}$. So $\mathcal{A} \subseteq \mathcal{L}_{\bf c}^{(r)}$ and $\mathcal{B} \subseteq \mathcal{L}_{\bf d}^{(s)}$. As is explained in Remark~\ref{seqlab}, $\mathcal{A}$ and $\mathcal{B}$ are cross-intersecting, $|\mathcal{L}_{\bf c}^{(r)}| = |S_{\bf c}^{(r)}|$, $|\mathcal{L}_{\bf d}^{(s)}| = |S_{\bf d}^{(s)}|$, $|\mathcal{A}| = |A|$ and $|\mathcal{B}| = |B|$. Thus, the result follows if we show that $|\mathcal{A}||\mathcal{B}| \leq |\mathcal{X}||\mathcal{Y}|$, and that equality holds only if and only if, for some $p \in \{h \in [\min\{m,n\}] \colon c_h = c_1, d_h = d_1\}$ and some $q \in [c_p]$, $\mathcal{A} = \{A \in \mathcal{L}_{\bf c}^{(r)} \colon (p,q) \in A\}$ and $B = \{B \in \mathcal{L}_{\bf d}^{(s)} \colon (p,q) \in B\}$; the sufficiency condition is trivial.

Let $\mathcal{G} = {[m] \choose \leq r}$. Let $v \colon \mathcal{G} \rightarrow \mathbb{N}$ such that for each $G \in \mathcal{G}$,
\[v(G) = \left| \left\{ L \in \mathcal{L}_{\bf c}^{(r)} \colon L \cap ([m] \times [1]) = G \times [1] \right\} \right|.\]
Let $\mathcal{H} = {[n] \choose \leq s}$. Let $w \colon \mathcal{H} \rightarrow \mathbb{N}$ such that for each $H \in \mathcal{H}$,
\[w(H) = \left| \left\{ L \in \mathcal{L}_{\bf d}^{(s)} \colon L \cap ([n] \times [1]) = H \times [1] \right\} \right|.\]
Let $l = \max\{m,n\}$ and $h = \max\{c_m, d_n\}$. Let
\begin{align} \mathcal{A}^* &= \Gamma_{l,h} \circ \dots \circ \Gamma_{l,2} \circ \dots \circ \Gamma_{2,h} \circ \dots \circ \Gamma_{2,2} \circ \Gamma_{1,h} \circ \dots \circ \Gamma_{1,2}(\mathcal{A}),\nonumber \\
\mathcal{B}^* &= \Gamma_{l,h} \circ \dots \circ \Gamma_{l,2} \circ \dots \circ \Gamma_{2,h} \circ \dots \circ \Gamma_{2,2} \circ \Gamma_{1,h} \circ \dots \circ \Gamma_{1,2}(\mathcal{B}). \nonumber
\end{align}
Now let
%
\begin{align} \mathcal{C} &= \left\{G \in \mathcal{G} \colon E \cap ([m] \times [1]) = G \times [1] \mbox{ for some } E \in \mathcal{A}^* \right\}, \nonumber \\
\mathcal{D} &= \left\{H \in \mathcal{H} \colon F \cap ([n] \times [1]) = H \times [1] \mbox{ for some } F \in \mathcal{B}^* \right\}. \nonumber
\end{align}
So $\mathcal{C} \subseteq \mathcal{G}$, $\mathcal{D} \subseteq \mathcal{H}$, and by Corollary~\ref{deltacor}, $\mathcal{C}$ and $\mathcal{D}$ are cross-intersecting. We have $\mathcal{A}^* \subseteq \bigcup_{C \in \mathcal{C}} \{L \in \mathcal{L}_{\bf c}^{(r)} \colon L \cap ([m] \times [1]) = C \times [1]\}$ and $\mathcal{B}^* \subseteq \bigcup_{D \in \mathcal{D}} \{L \in \mathcal{L}_{\bf d}^{(s)} \colon L \cap ([n] \times [1]) = D \times [1]\}$. So 
\begin{equation} |\mathcal{A}^*| \leq \sum_{C \in \mathcal{C}} v(C) = v^{(\mathcal{G})}(\mathcal{C}) \quad \mbox{and} \quad |\mathcal{B}^*| \leq \sum_{D \in \mathcal{D}} w(D) = w^{(\mathcal{H})}(\mathcal{D}). \label{maingen.1}
\end{equation} 
Since $|\mathcal{A}| = |\mathcal{A}^*|$ and $|\mathcal{B}| = |\mathcal{B}^*|$, we therefore have 
\begin{equation} |\mathcal{A}| \leq v^{(\mathcal{G})}(\mathcal{C}) \quad \mbox{and} \quad |\mathcal{B}| \leq w^{(\mathcal{H})}(\mathcal{D}). \label{maingen.15}
\end{equation}
Let $\mathcal{I} = \{G \in \mathcal{G} \colon 1 \in G\}$ and $\mathcal{J} = \{H \in \mathcal{H} \colon 1 \in H\}$. By Lemma~\ref{lemmaweighted} and Theorem~\ref{xintweight}, 
\begin{equation} v^{(\mathcal{G})}(\mathcal{C}) w^{(\mathcal{H})}(\mathcal{D}) \leq v^{(\mathcal{G})}(\mathcal{I}) w^{(\mathcal{H})}(\mathcal{J}). \label{maingen.2}
\end{equation}  
Now
\begin{align} v^{(\mathcal{G})}(\mathcal{I}) &= \left( \sum_{I \in \mathcal{I}} v(I) \right) = \left(\sum_{I \in \mathcal{I}} \left| \left\{L \in \mathcal{L}_{\bf c}^{(r)} \colon L \cap ([m] \times [1]) = I \times [1] \right\} \right| \right) \nonumber \\
&= \left|\bigcup_{I \in \mathcal{I}} \left \{L \in \mathcal{L}_{\bf c}^{(r)} \colon L \cap ([m] \times [1]) = I \times [1] \right\} \right| = |\mathcal{X}| \nonumber 
\end{align}
and similarly $w^{(\mathcal{H})}(\mathcal{J}) = |\mathcal{Y}|$. Together with (\ref{maingen.15}) and (\ref{maingen.2}), this gives us $|\mathcal{A}||\mathcal{B}| \leq |\mathcal{X}||\mathcal{Y}|$. Suppose equality holds. Then all the inequalities in (\ref{maingen.1})--(\ref{maingen.2}) are equalities. The equalities in (\ref{maingen.1}) imply that $\mathcal{A}^* = \bigcup_{C \in \mathcal{C}} \{L \in \mathcal{L}_{\bf c}^{(r)} \colon L \cap ([m] \times [1]) = C \times [1]\}$ and $\mathcal{B}^* = \bigcup_{D \in \mathcal{D}} \{L \in \mathcal{L}_{\bf d}^{(s)} \colon L \cap ([n] \times [1]) = D \times [1]\}$. By Theorem~\ref{xintweight}, equality in (\ref{maingen.2}) gives us that for some $p \in [m] \cap [n]$, $\mathcal{C} = \mathcal{G}(p)$ and $\mathcal{D} = \mathcal{H}(p)$. 
It follows that $\mathcal{A}^* = \{L \in \mathcal{L}_{\bf c}^{(r)} \colon (p,1) \in L\}$ and $\mathcal{B}^* = \{L \in \mathcal{L}_{\bf d}^{(s)} \colon (p,1) \in L\}$. By Lemma~\ref{ss_cross_int}, $\mathcal{A} = \{L \in \mathcal{L}_{\bf c}^{(r)} \colon (p,q) \in L\}$ and $\mathcal{B} = \{L \in \mathcal{L}_{\bf d}^{(s)} \colon (p,q) \in L\}$ for some $q \in [c_p] \cap [d_p]$. So $\mathcal{A}$ is a star of $\mathcal{L}_{\bf c}^{(r)}$ with centre $(p,q)$, and $\mathcal{B}$ is a star of $\mathcal{L}_{\bf d}^{(s)}$ with centre $(p,q)$. Now clearly $\mathcal{X}$ is a star of $\mathcal{L}_{\bf c}^{(r)}$ of maximum size, and $\mathcal{Y}$ is a star of $\mathcal{L}_{\bf d}^{(s)}$ of maximum size. Thus, since $|\mathcal{A}||\mathcal{B}| = |\mathcal{X}||\mathcal{Y}|$, $|\mathcal{A}| = |\mathcal{X}|$ and $|\mathcal{B}| = |\mathcal{Y}|$. So $\mathcal{A}$ is a star of $\mathcal{L}_{\bf c}^{(r)}$ of maximum size, and hence we must have $c_p = c_1$. Similarly, $d_p = d_1$.~\hfill{$\Box$}\\
\\
\textbf{Proof of Theorem~\ref{main}.} Since $|A| \leq |S_{\bf c}|$ and $|B| \leq |S_{\bf c}|$, the result is trivial if $c_1 = 1$.

If $c_1 \geq 3$, then the result is given by Theorem~\ref{mainpartial} with $r = s = m = n$ and ${\bf c} = {\bf d}$.

Finally, suppose $c_1 = 2$. Let
\begin{align} \mathcal{A} &= \{\{(1,a_1), \dots, (n,a_n)\} \colon (a_1, \dots, a_n) \in A\}, \nonumber \\
\mathcal{B} &= \{\{(1,b_1), \dots, (n,b_n)\} \colon (b_1, \dots, b_n) \in B\}. \nonumber
\end{align}
So $\mathcal{A}, \mathcal{B} \subseteq \mathcal{L}_{\bf c}$, $|\mathcal{A}| = |A|$, $|\mathcal{B}| = |B|$, and as is explained in Remark~\ref{seqlab}, $\mathcal{A}$ and $\mathcal{B}$ are cross-intersecting. Thus, since $|\mathcal{L}_{\bf c}| = |S_{\bf c}|$, the result follows if we show that $|\mathcal{A}||\mathcal{B}| \leq \left( \frac{1}{c_1}|\mathcal{L}_{\bf c}| \right)^2$. Let mod$^*$ be as in the proof of Theorem~\ref{mainpartialgen}. Let $\theta : \mathcal{L}_{\bf c} \rightarrow \mathcal{L}_{\bf c}$ such that $\theta(E) = \{(i,(j+1) \mbox{ mod$^*$ } c_i) \colon (i,j) \in E\}$ for each $E \in \mathcal{L}_{\bf c}$. Clearly, $\theta$ is a bijection, and $\theta(E) \cap E = \emptyset$ for each $E \in \mathcal{L}_{\bf c}$. Thus, since $\mathcal{A}$ and $\mathcal{B}$ are cross-intersecting, $\theta(C) \notin \mathcal{B}$ for each $C \in \mathcal{A}$, and hence $|\mathcal{B}| \leq |\mathcal{L}_{\bf c}| - |\mathcal{A}|$. Since $0 \leq \left(|\mathcal{A}| - \frac{1}{2}|\mathcal{L}_{\bf c}| \right)^2 = |\mathcal{A}|^2 - |\mathcal{A}||\mathcal{L}_{\bf c}| + \frac{1}{4}|\mathcal{L}_{\bf c}|^2$, we have $|\mathcal{A}|(|\mathcal{L}_{\bf c}| - |\mathcal{A}|) \leq \left( \frac{1}{2}|\mathcal{L}_{\bf c}| \right)^2$ and hence $|\mathcal{A}||\mathcal{B}| \leq \left( \frac{1}{c_1}|\mathcal{L}_{\bf c}| \right)^2$.~\hfill{$\Box$}

\end{document}